\documentclass[reqno]{article}
\usepackage{amsmath, amsthm, graphicx}
\theoremstyle{plain}

\usepackage{setspace}
\usepackage{ifthen}

\begin{document}
\bibliographystyle{amsalpha}
\title{Another Riemann-Farey Computation}
\author{Scott B. Guthery\\sguthery@mobile-mind.com}
\maketitle

The Riemann hypothesis is true if and only if
\begin{equation}\label{EQ1}
R(m) = \sum_{i=2}^{T_{m}}\left({F_m(i) - \frac{i}{n}}\right)^2 =
O(m^{-1+\epsilon})
\end{equation}
where $F_m(i)$ is the $i^{th}$ element in the Farey sequence of
order $m$ and
$$
T_m=\sum_{k=2}^m \phi(k)\text{.}
$$

Let $P_m(k)$ be sum of the $\phi(k)$ terms in (\ref{EQ1}) with Farey
denominator $k$ so that

\begin{equation}\label{EQ2}
R(m) = \sum_{k=2}^m P_m(k)
\end{equation}

\pagebreak

\begin{figure}[h!]
\begin{center}
\leavevmode
  \includegraphics[width=100mm]{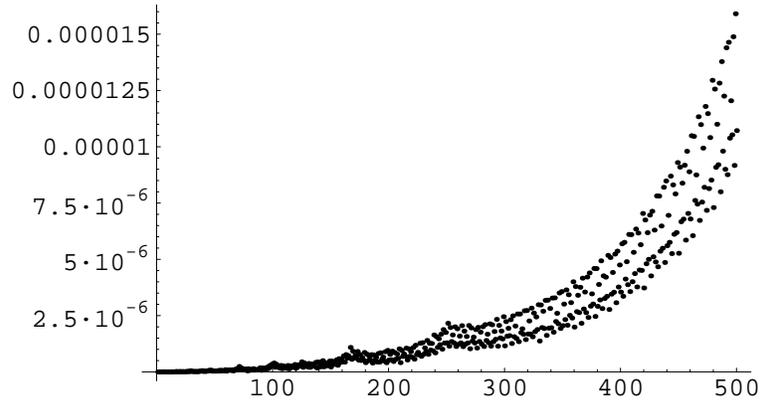}\\
  \end{center}\caption{$P_{m}(k)$ for $m=500$}\label{FIG1}
\end{figure}

Figure~\ref{FIG1} is a plot of $P_m(k)$ for $m = 500$.

Let $Q_{m}(k,i)$ denote the term with numerator $i$ in $P_{m}(k)$ so
that

\begin{equation}\label{EQ3}
P_m(k) = \sum_{i=1}^m Q_{m}(k,i)
\end{equation}

and we note in passing that for all but small values of $m$ and $k$,

\begin{equation}\label{EQ4}
Q_{m}(k,1) \gg Q_{m}(k,j)
\end{equation}

for $j>1$.

\pagebreak

\begin{figure}[h!]
\begin{center}
\leavevmode
  \includegraphics[width=100mm]{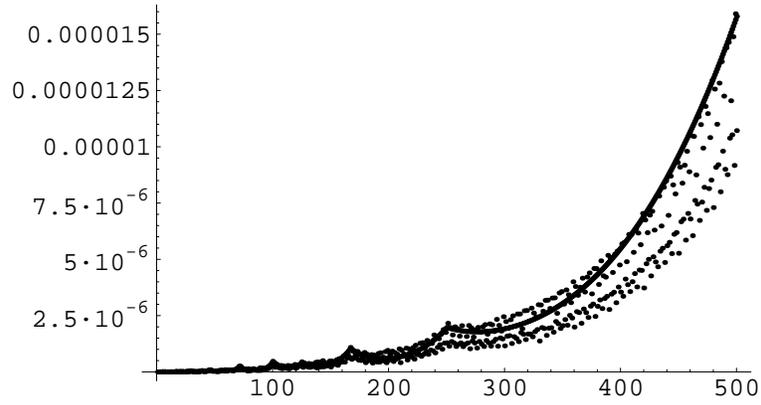}\\
  \end{center}\caption{$P_{500}(k)$ and $4 Q_{500}(k,1)$}\label{FIG3}
\end{figure}

The solid line in Figure~\ref{FIG3} is four times $Q_{500}(k,1)$
plotted with $P_{500}(k)$.

Thus, we are led to consider

\begin{equation}\label{EQ5}
R(m) < C \sum_{k=2}^m Q_m(k,1)
\end{equation}

\pagebreak

For $83 \leq k \leq 500$, $Q_{500}(k,1)$ is given by

\begin{equation}\label{EQ6}
Q_{500}(k,1) =
\begin{cases}

\left(\frac{501-k}{76116} - \frac{1}{k}\right)^2,&

250 \leq k \leq 500\\

\left(\frac{249+2(251-k)}{76116} - \frac{1}{k}\right)^2,&

167 \leq k < 250\\

\left(\frac{413+4(168-k)}{76116} - \frac{1}{k}\right)^2,&

125 \leq k < 167\\

\left(\frac{579+6(126-k)}{76116} - \frac{1}{k}\right)^2,&

100 \leq k < 125\\

\left(\frac{725+10(101-k)}{76116} - \frac{1}{k}\right)^2,&

83 \leq k < 100
\end{cases}
\end{equation}

Each case in (\ref{EQ5}) corresponds to a constant step size between
adjacent Farey numbers of the form $\frac{1}{k}$. This step size is
the factor of $(a-k)$ in the case.

\pagebreak

Considering the top case in equation (\ref{EQ4}), for $\frac{m}{2}
\leq k \leq m$ we have in general

\begin{equation}
Q_m(k,1) = \left(\frac{m+1-k}{\sum_{j=2}^m
\phi(j)}-\frac{1}{k}\right)^2
\end{equation}

Substituting the approximation $\frac{3 m^2}{\pi ^2}$ for the
totient sum and integrating from $m/2$ to $m$, we have

\begin{equation}
\begin{split}
\hat{R}(m) = \frac{12\pi^4 + 6m(\pi^4-24\pi^2 \log{2})+
m^2(216+\pi^4-72\pi^2(\log{4}-1))}{216 m^3}\\
\approx \frac{0.180 m^2-1.855 m+5.41}{m^3}
\end{split}
\end{equation}

with

\begin{equation}
\lim_{m\to \infty} \hat{R}(m)/m^{-1+\epsilon} = 0
\end{equation}

 \pagebreak

In general, each of the $n$ terms in (\ref{EQ1}) can be represented
as

\begin{equation}\label{EQ7}
\left(\frac{a + b k}{\sum_{j=2}^m \phi(j)} - \frac{i}{k}\right)^2
\end{equation}

for $[m/c] \le k < [m/c]+1$ and $a$, $b$ and $c$ depending on $m$,
$k$ and $i$.

As above, substituting the approximation $\frac{3 m^2}{\pi ^2}$ for
the totient sum, integrating from $\frac{m}{c}$ to $\frac{m}{c}+1$
for each term and then summing over the $n$ terms, we have

\begin{equation}
\begin{split}
R(m) < \hat{R}(m) = C \sum_{k,i}{\frac{27 c^4 i^2 m^3 - 18 b c^2 i
m^2 (c+m) \pi^2 +
18 a c^2 i m^2(c+m)\pi^2\log{\frac{m}{(m + c)}}}{27 c^2(m + c) m^4}}\\
+\frac{(m + c)(3 a^2 c^2 + 3a b c (c+2m) + b^2(c^2+3 c m + 3 m^2))
\pi^4}{27 c^2(m + c) m^4}
\end{split}
\end{equation}

where the constant $C$ accounts for the totient sum approximation
and

\begin{equation}
\lim_{m\to \infty}\hat{R}(m)/m^{-1+\epsilon} = 0
\end{equation}

\end{document}